\documentclass{article}

\usepackage{amsmath}
\usepackage{amsfonts}
\usepackage{amssymb}
\usepackage{graphicx}
\usepackage{color}
\usepackage{url}
\usepackage[latin1]{inputenc}
\usepackage[T1]{fontenc}

\newtheorem{theorem}{Theorem}[section]

\newtheorem{definition}{Definition}[section]

\newtheorem{lemma}{Lemma}[section]

\newtheorem{remarks}{Remarks}[section]

\begin{document}

\title{The Navier-Stokes problem modified by an absorption term}

\author{Hermenegildo Borges de Oliveira\\
\bigskip
holivei@ualg.pt\\
FCT - Universidade do Algarve\\
Campus de Gambelas\\
8005-139 Faro\\
Portugal}

\date{}

\maketitle

\begin{abstract}
In this work we consider the Navier-Stokes problem modified by the absorption term $|\textbf{u}|^{\sigma-2}\textbf{u}$, where $\sigma>1$, which is  introduced in the momentum equation. %
For this new problem, we prove the existence of weak solutions for any dimension $N\geq 2$ and its uniqueness
for $N=2$. %
Then we prove that, for zero body forces, the weak solutions extinct in a finite time if $1<\sigma<2$, exponentially decay in time if $\sigma=2$ and decay with a power-time rate if $\sigma>2$. %
We prove also that for a general non-zero body forces, the weak solutions exponentially decay in time for any $\sigma>1$.
In the special case of a suitable forces field which vanishes at some instant, we prove that the weak solutions extinct at the same instant provided $1<\sigma<2$.

\end{abstract}


\bigskip


\section{Introduction}

\subsection{Basic equations}

From the basic principles of Fluid Mechanics, it is well known that, in isothermal motions of incompressible fluids, the velocity field and pressure are determined from:
\begin{itemize}
\item the incompressibility condition
\begin{equation}\label{inc-cond}
\mathrm{div}\textbf{u}=0;
\end{equation}
\item the conservation of mass
\begin{equation}\label{cons-mass}
\frac{\partial\,\rho}{\partial\,t}+\mathrm{div}(\rho\textbf{u})=0;
\end{equation}
\item the conservation of momentum
\begin{equation}\label{cons-mom}
\rho\left(\frac{\partial\,\textbf{u}}{\partial\,t}+(\textbf{u$\cdot\nabla$})\textbf{u}\right)=\rho\,\textbf{f}+\mathrm{div}\textbf{S}.
\end{equation}
\end{itemize}
The notation used in (\ref{inc-cond})-(\ref{cons-mom}) is well
known: $\textbf{u}$ is the velocity field, $p$ is the pressure,
$\rho$ is the density and $\textbf{f}$ is the forcing term. %
For Newtonian fluids, the stress tensor $\textbf{S}$ is given by the Stokes law %
\begin{equation}\label{stress-t}
\textbf{S}=-p\textbf{I}+2\mu\textbf{D}\,,\quad
\textbf{D}=\frac{1}{2}\left(\textbf{$\nabla$\,u}+\textbf{$\nabla$\,u}^T\right)\,,
\end{equation}
where $\mu$ is the  dynamical (constant) viscosity and $\textbf{D}$ is the strain tensor. For homogeneous fluids, the density is regarded as
constant. Therefore, we can replace the continuity equation (\ref{cons-mass}) by
its incompressible form (\ref{inc-cond}). In consequence, when these simplifying
features are present, the equations for a linearly viscous, homogeneous, incompressible fluid reduce to the
following system
\begin{equation}\label{NSE}
\frac{\partial\,\textbf{u}}{\partial\,t}+(\textbf{u$\cdot\nabla$})\textbf{u}=\textbf{f}-\frac{1}{\rho}\textbf{$\nabla$}p+\nu\triangle\textbf{u} \end{equation}
\begin{equation}\label{inc}
\mathrm{div}\textbf{u}=0\,.
\end{equation}
Here $\nu=\mu/\rho$ is a constant corresponding to the denominated
kinematics viscosity and, because the fluid is homogeneous, $\rho$ is a constant corresponding to the initial density.
In the literature, the system of equations (\ref{NSE})-(\ref{inc}) is known  as the incompressible Navier-Stokes equations, or only the Navier-Stokes equations (NSE). System (\ref{NSE})-(\ref{inc}) must be supplemented with boundary conditions characterizing the flow on the boundary of the domain occupied by the fluid and by initial conditions determining the initial state of the flow at the beginning of the time interval. The question of initial conditions is immediately understood from the physical point of view, but as for the boundary conditions is much more delicate and would require a detailed discussion. See Serrin~\cite{Serrin-1959} for a detailed derivation of the Navier-Stokes equations and for the explanation of the initial and possible different boundary conditions.

\subsection{The modified Navier-Stokes problem}\label{sect-TMNS}

Let us consider a general cylinder
$$Q_T:=\Omega\times(0,T)\subset\mathbb{R}^N\times\mathbb{R}^+\,,\quad\mbox{with}\quad
\Gamma_T:\partial\Omega\times(0,T)\,,$$ %
where $\Omega$ is a bounded domain with a compact
boundary $\partial\Omega$. Although the dimensions of physical interest are $N=2$ and $N=3$, here we shall consider a general dimension $N\geq 2$. %
In this work we consider the following modified NSE
\begin{equation}\label{NSE-2}
\frac{\partial\,\textbf{u}}{\partial\,t}+(\textbf{u$\cdot\nabla$})\textbf{u}=\textbf{f}
-\textbf{$\nabla$}p+\nu\triangle\textbf{u}-\alpha|\textbf{u}|^{\sigma-2}\textbf{u}
\quad\mbox{in}\quad
Q_T\,,
\end{equation}
\begin{equation}\label{inc-2}
\mathrm{div}\textbf{u}=0\quad\mbox{in}\quad Q_T
\end{equation}
supplemented by the initial and boundary conditions
\begin{equation}\label{i-c}
\textbf{u}=\textbf{u}_0\quad\mbox{when}\quad t=0\quad \mbox{in}\quad\Omega
\end{equation}
\begin{equation}\label{b-c}
\textbf{u}=\textbf{0}\quad\mbox{on}\quad \Gamma_T\,.
\end{equation}
In (\ref{NSE-2}), $\alpha$ and $\sigma$ are positive constants with $\sigma>1$. %
Notice that, for the sake of simplicity,  we have assumed $\rho\equiv
1$ in (\ref{NSE-2}). %
In the sequel (\ref{NSE-2})-(\ref{b-c}) shall be denominated as the modified NS problem. %
The motivation for the consideration of the absorption term $|\textbf{u}|^{\sigma-2}\textbf{u}$ in (\ref{NSE-2}) is
purely mathematical and goes back to the works of Benilan \emph{et
al.} \cite{BBC-1975}, D\'iaz and Herrero \cite{DH-1981}, and Bernis
\cite{Bernis-1984,Bernis-1986}. %
There, was studied the importance of a similar absorption term to prove
qualitative properties related with compact supported solutions, or
solutions which exhibit finite speed of propagations, or which
extinct in time, for different initial and or boundary value
problems. %
A possible physical justification for the absorption term in (\ref{NSE-2}), is the consideration, in the
momentum equation, of a forces field
\begin{equation}\label{force-feedback}
\textbf{h}(\textbf{u})=\textbf{f}-\alpha|\textbf{u}|^{\sigma-2}\textbf{u}\,,
\end{equation}
where $\textbf{f}$ is a given forces field. %
Notice that such forces field depends, in a sublinear way, on the own velocity $\textbf{u}$.
In a certain sense, such forces field may be considered, from the physical point of view, as a feedback field. Therefore,
(\ref{NSE-2}) can be considered as the equation of motion for a
certain fluid. In this way, the purpose of this work is to study the
response of such a fluid motion undergone to a forces field
satisfying (\ref{force-feedback}). Specifically, we want to know if
such a forces field is responsible for stopping the fluid, driven by
(\ref{NSE-2})-(\ref{b-c}). This issue addresses us for the important question about
the decays of the solutions of the Navier-Stokes equations.
With this respect, it should be remarked that these questions have been studied by many authors, among many others, by
Horgan and Wheeler~\cite{HW-1978}, Wiegner~\cite{W-1987}, Ames and Payne~\cite{AP-1989}, Schonbek~\cite{Schonbek-1991}, Borchers and Miyakawa~\cite{BM-1992}, Kozono and Ogawa~\cite{KO-1993}, Takahashi~\cite{Taka-1999}, Enomoto and Shibata~\cite{ES-2005} and Bae and Jin~\cite{BJ-2005}. In the references~\cite{HW-1978,AP-1989} are established the exponential spacial decay for the solutions of the stationary Navier-Stokes equations. The other references \cite{W-1987,Schonbek-1991,BM-1992,KO-1993,Taka-1999,ES-2005,BJ-2005} are concerned with space and time power decays in different norms for the solutions and its derivatives of the Navier-Stokes equations supplemented with suitable conditions on the initial data or in the forces field. In spite of many work in this field, so far, and to the best of our knowledge, there are no results establishing the extinction of the solutions of the Navier-Stokes equations in a finite time or in space.

\subsection{Mathematical framework}\label{frame}

\textbf{Notation.} The notation used throughout this text is largely
standard in Mathematical Analysis and in particular in Mathematical Fluid Mechanics - see, \emph{e.g.}, Lions~\cite{Lions-1969}, Temam~\cite{Temam-1979} or Galdi~\cite{Galdi-1994,Galdi-2000}. We distinguish vectors from scalars by
using boldface letters. For functions and function spaces we will
use this distinction as well. The symbol $C$ will denote a generic
constant - almost the time a positive constant, whose value will not
be specified; it can change from one inequality to another. The
dependence of $C$ on other constants or parameters will always be
clear from the exposition. %
Sometimes we will use subscripted letters attached to $C$ to relate a constant with the result it is derived from. %
In this
article, the notation $\Omega$ stands always for a domain,
\emph{i.e.}, a connected open subset of $\mathbb{R}^N$, whose
compact boundary is denoted by $\partial\Omega$. %

\bigskip\noindent\textbf{Function spaces.} %
Let $1\leq p\leq \infty$ and $\Omega\subset\mathbb{R}^N$ be a domain. %
We shall use the classical Lebesgue spaces $\mathrm{L}^p(\Omega)$, whose norm is
denoted by $\|\cdot\|_{\mathrm{L}^p(\Omega)}$. %
For any nonnegative $k$,
$\mathrm{W}^{k,p}(\Omega)$ denotes the Sobolev space of all
functions $u\in\mathrm{L}^p(\Omega)$ such
that the weak derivatives $\mathrm{D}^{\alpha}u$ exist, in the generalized sense, and are in
$\mathrm{L}^p(\Omega)$ for any multi-index $\alpha$ such that
$0\leq |\alpha|\leq k$. The norm in $\mathrm{W}^{k,p}(\Omega)$ is denoted by
$\|\cdot\|_{\mathrm{W}^{k,p}(\Omega)}$.
The associated trace spaces are denoted by
$\mathrm{W}^{k-1/p,p}(\partial\Omega)$.
Given $T>0$ and a Banach space $X$,
$\mathrm{L}^p(0,T;X)$ and $\mathrm{W}^{k,p}(0,T;X)$ denote
the usual Lebesgue and Sobolev spaces used in evolutive problems,
with norms denoted by $\|\cdot\|_{\mathrm{L}^p(0,T;X)}$ and
$\|\cdot\|_{\mathrm{W}^{k,p}(0,T;X)}$. %
The corresponding spaces of vector-valued
functions are denoted by boldface letters. All these spaces are
Banach spaces and the Hilbert framework corresponds to $p=2$.
In the last case, we use the abbreviations $\mathrm{W}^{k,2}=\mathrm{H}^{k}$ and
$\mathrm{W}^{k-1/2,2}=\mathrm{H}^{k-1/2}$.

\bigskip\noindent\textbf{Auxiliary results.}
Throughout this text we will make reference, more than once, to the
following inequalities:\\
(1) Algebraic inequality - for every $\alpha$, $\beta\in\mathbb{R}$
and every $A$, $B\geq 0$,
\begin{equation}\label{algebraic}
A^{\alpha}B^{\beta}\leq (A+B)^{\alpha+\beta};
\end{equation}
(2) Young's inequality - for every $a,\ b\geq 0$, $\varepsilon>0$
and $1<p,\ q<\infty$ such that $1/p+1/q=1$,
$$a\,b\leq\varepsilon a^p+C(\varepsilon)b^q\,,\quad C(\varepsilon)=(\varepsilon\,p)^{-q/p}q^{-1}\,.$$
If $p=q=2$, this is known as Cauchy's inequality.\\
(3) Hölder's inequality - for every $u\in\mathrm{L}^p(\Omega)$, $v\in
\mathrm{L}^q(\Omega)$, with $1\leq p,\ q\leq \infty$ such that $1/p+1/q=1$,
$$\int_{\Omega}u\,v\,d\textbf{x}\leq
\|u\|_{\mathrm{L}^p(\Omega)}\|v\|_{\mathrm{L}^q(\Omega)}.$$ %

An important result that will be used in the sequel is the famous Gagliardo-Nirenberg-Sobolev inequality.

\begin{lemma}\label{L-Nirenberg}
Let $\Omega$ be a domain of $\mathbb{R}^N$, $N\geq 1$, with a
compact boundary $\partial\Omega$. Assume that $u\in
\mathrm{W}^{1,p}_0(\Omega)$. Then, for every fixed number $r\geq 1$ there
exists a constant $C$ depending only on $N$, $p$, $r$ such that
\begin{equation}\label{Nirenberg}
\left\|u\right\|_{\mathrm{L}^q(\Omega)}\leq C
\|\mathbf{\nabla}u\|_{\mathrm{L}^p(\Omega)}^{\theta}\|u\|_{\mathrm{L}^r(\Omega)}^{1-\theta},
\end{equation}
where $p,\ q\geq 1$, are linked by
\begin{equation}\label{theta}
\theta=\left(\frac{1}{r}-\frac{1}{q}\right)\left(\frac{1}{N}-\frac{1}{p}+\frac{1}{r}\right)^{-1}\,,
\end{equation}
and their admissible range is:
\begin{enumerate}
\item If $N=1$, $q\in[r,\infty]$, $\theta\in\left[0,\frac{p}{p+r(p-1)}\right]$, $C=[1+(p-1)/pr]^{\theta}$; %
\item If $p<N$, $q\in\left[\frac{Np}{N-p},r\right]$ if $r\geq \frac{Np}{N-p}$ and %
$q\in\left[r,\frac{Np}{N-p}\right]$ if $r\leq \frac{Np}{N-p}$,  $\theta\in[0,1]$ and $C=[(N-1)p/(N-p)]^{\theta}$; \\ %
\item If $p\geq N>1$, $q\in[r,\infty)$, $\theta\in\left[0,\frac{Np}{Np+r(p-N)}\right)$ and $C=\max\{q(N-1)/N,1+(p-1)pr\}^{\theta}$.
\end{enumerate}
\end{lemma}
When $\theta=1$, (\ref{Nirenberg}) is known as the Sobolev inequality and, in this case, if $q=p=2$, then (\ref{Nirenberg}) is usually denominated as the Poincaré inequality. %
This result is valid whether the domain $\Omega$ is bounded or not and notice the constant $C$ does not depend on $\Omega$. See the proof in Ladyzhenskaya~\emph{et al.}~\cite[p. 62]{LSU-1967}. %
The following generalization of (\ref{Nirenberg}) to higher-order derivatives is also possible (see Nirenberg~\cite{Nirenberg-1959}),
\begin{equation}\label{Nirenberg-ho}
\left\|\mathrm{D}^ju\right\|_{\mathrm{L}^q(\Omega)}\leq C
\|\mathrm{D}^ku\|_{\mathrm{L}^p(\Omega)}^{\theta}\|u\|_{\mathrm{L}^r(\Omega)}^{1-\theta}\,,
\end{equation}
where additionally  $0\leq j<k$ and $k\geq 1$, and now %
$$\theta=\left(\frac{1}{r}+\frac{j}{N}-\frac{1}{q}\right)\left(\frac{k}{N}-\frac{1}{p}+\frac{1}{r}\right)^{-1}\,.$$

The results written above can be easily generalized for vector-valued functions. Moreover, for time-dependent functions $u(x,t)$, they still hold for a.a. $t\in[0,T]$.

Other very important auxiliary result which will be of the utmost importance to handle the absorption term is written in the following lemma.

\begin{lemma}
For all $p\in(1,\infty)$ and $\delta\geq 0$, there exist constants $C_1$ and $C_2$, depending on $p$ and $N$,
such that for all $\textbf{$\xi$},\ \textbf{$\eta$}\in\mathbb{R}^N$, $N\geq 1$,
\begin{equation}\label{lin-xi-eta}
\left||\textbf{$\xi$}|^{p-2}\textbf{$\xi$}-|\textbf{$\eta$}|^{p-2}\textbf{$\eta$}\right|\leq
C_1|\textbf{$\xi$}-\textbf{$\eta$}|^{1-\delta}\left(|\textbf{$\xi$}|+|\textbf{$\eta$}|\right)^{p-2+\delta}
\end{equation}
and
\begin{equation}\label{gin-xi-eta}
\left(|\textbf{$\xi$}|^{p-2}\textbf{$\xi$}-|\textbf{$\eta$}|^{p-2}\textbf{$\eta$}\right)\cdot(\textbf{$\xi$}-\textbf{$\eta$})\geq
C_2|\textbf{$\xi$}-\textbf{$\eta$}|^{2+\delta}\left(|\textbf{$\xi$}|+|\textbf{$\eta$}|\right)^{p-2-\delta}\,.
\end{equation}
\end{lemma}
See Barret and Liu~\cite{BL-1994} for the proof and also the references cited therein for other forms of (\ref{lin-xi-eta}) and (\ref{gin-xi-eta}).

\section{Weak formulation}

In this section, we establish existence and uniqueness results for
the modified NS problem (\ref{NSE-2})-(\ref{b-c}).  Notice that, in the limit case $\alpha=0$ in
(\ref{NSE-2}), we fall in the classical NS system (\ref{NSE})-(\ref{inc}).
In this case, it is well known that the corresponding problem has a
weak solution which is unique only if $N=2$
(see \emph{e.g.} Lions~\cite{Lions-1969} and Galdi~\cite{Galdi-2000}). %
In order to define the notion of a weak solution to the
problem (\ref{NSE-2})-(\ref{b-c}), let us introduce the function spaces largely used
in Mathematical Fluid Mechanics:

\begin{equation}\label{space-V-cal}
\mathcal{V}:=\{\textbf{v}\in\textbf{C}_0^{\infty}(\Omega):\mathrm{div}\textbf{v}=0\}\,;
\end{equation}
\begin{equation}\label{space-H}
\textbf{H}:=\mbox{closure of $\mathcal{V}$ in $\textbf{L}^2(\Omega)$}\,;
\end{equation}
\begin{equation}\label{space-V-s}
\textbf{V}_{s}:=\mbox{closure of $\mathcal{V}$ in $\textbf{H}^s(\Omega)$},\quad s\geq 1\,;
\end{equation}
where, for simplicity, we can assume $s$ as the smaller integer not lesser than $N/2$, to avoid the complicated Sobolev spaces
with $s$ non-integer. %
$\textbf{H}$ is endowed with the $\textbf{L}^2(\Omega)$ inner product and norm, and, for any positive integer $s$, $\textbf{V}_s$ is endowed with the inner product
$$(\textbf{u},\textbf{v})_{\textbf{V}_s}=
\sum_{|\alpha|=s}\int_{\Omega}\textbf{D}^{\alpha}\textbf{u}\cdot\textbf{D}^{\alpha}\textbf{v}\,d\textbf{x}$$
and with the associated norm. From (\ref{Nirenberg}), we see
that this norm is equivalent to the $\textbf{H}_0^s(\Omega)$ norm and, in consequence, we have the (isometrically isomorphic) identification
$\textbf{V}'_s=\textbf{H}^{-s}(\Omega)$. %
$\textbf{H}$ and $\textbf{V}_{s}$ are Hilbert spaces and, for $s=1$, which happens when the dimension is $N=2$, we simply denote, as usual,
$\textbf{V}_1$ by $\textbf{V}$. In this case, the $\textbf{V}$ inner product reads in the usual form
$$(\textbf{u},\textbf{v})_{\textbf{V}}=\int_{\Omega}\mathbf{\nabla\,u:\nabla\,v}\,d\textbf{x}\,.$$ %
From the theory of distributions, we know that
\begin{equation}\label{distr-incl}
\textbf{V}_{s}\hookrightarrow\textbf{V}\hookrightarrow\textbf{H}=\textbf{H}'\hookrightarrow\textbf{V}'\hookrightarrow\textbf{V}_{s}'\,,\quad s>1\,.
\end{equation}
Moreover the compact imbedding $\textbf{H}^1_0(\Omega)\hookrightarrow\textbf{L}^2(\Omega)$ implies that the imbedding $\textbf{V}_{s}\hookrightarrow\textbf{H}$ is also compact for any $s\geq 1$ (see Lions~\cite[pp. 67-77]{Lions-1969}). %
For the theory of these function spaces, see Lions~\cite{Lions-1969}, Temam~\cite{Temam-1979} and Galdi~\cite{Galdi-1994}. %
The notion of weak solution for the problem (\ref{NSE-2})-(\ref{b-c}) follows in a standard manner.

\begin{definition}\label{def-w-sol-tp}
Let $\Omega$ be a bounded domain in $\textbf{R}^N$, $N\geq 2$. A vector field $\textbf{u}$
is a weak solution of the problem (\ref{NSE-2})-(\ref{b-c}), if $\sigma>1$ and:
\begin{enumerate}
\item $\textbf{u}\in\textrm{L}^2(0,T;\textbf{V})\cap\textrm{L}^{\sigma}(Q_T)\cap\textrm{L}^{\infty}(0,T;\textbf{H})$;
\item
   $\textbf{u}(\cdot,0)=\textbf{u}_0$ a.e. in $\Omega$;
\item For every $\textbf{v}\in\textbf{V}\cap\textbf{L}^N(\Omega)\cap\textbf{L}^{\sigma}(\Omega)$ and for a.a. $t\in(0,T)$,
\begin{equation}\label{w-s-s-h-t}
\begin{split}
\frac{d}{d\,t}\int_{\Omega}\textbf{u}(t)\cdot\textbf{v}\,d\,\textbf{x}+
\nu\int_{\Omega}\mathbf{\nabla}\textbf{u}(t):&\,\mathbf{\nabla}\textbf{v}\,d\,\textbf{x}+\int_{\Omega}\left[\left(\textbf{u}(t)\cdot\textbf{$\nabla$}\right)\textbf{u}(t)\right]\cdot\textbf{v}\,d\,\textbf{x} \\
+\,\alpha\int_{\Omega}|\textbf{u}(t)|^{\sigma-2}\textbf{u}(t)\cdot&\,\textbf{v}\,d\,\textbf{x}
=\int_{\Omega}\textbf{f}(t)\cdot\textbf{v}\,d\,\textbf{x}\,.
\end{split}
\end{equation}
\end{enumerate}
\end{definition}
Condition $\textbf{u}\in\textrm{L}^2(0,T;\textbf{V})$ expresses, in a certain sense, the incompressibility condition (\ref{inc-2}) and the boundary condition (\ref{b-c}), whereas condition $\textbf{u}\in\textrm{L}^{\infty}(0,T;\textbf{H})$, \emph{a priori}, does not seem to be strictly necessary or else to restrict the class of admissible weak solutions. Condition $\textbf{u}\in\textrm{L}^{\sigma}(Q_T)$ is a natural requirement to deal with the absorption term $|\textbf{u}|^{\sigma-2}\textbf{u}$. %
Condition $\textbf{u}(\cdot,0)=\textbf{u}_0$ a.e. in $\Omega$, should be interpreted in the sense that $\textbf{u}$ is $\textbf{L}^2(\Omega)$ weakly continuous for $t=0$, \emph{i.e.}
$$\lim_{t\to 0}\int_{\Omega}\left(\textbf{u}(t)-\textbf{u}_0\right)\cdot\textbf{v}\,d\,\textbf{x}=0\qquad \forall\ \textbf{v}\in\textbf{L}^{2}(\Omega)\,.$$
For $N$ and $\sigma\leq 4$, (\ref{w-s-s-h-t}) holds for every $\textbf{v}\in\textbf{V}$, because, due to Sobolev's inequality,
$\textbf{H}^1(\Omega)\hookrightarrow\textbf{L}^p(\Omega)$ for $p\leq 4$ and, in consequence, $\textbf{V}\cap\textbf{L}^p(\Omega)=\textbf{V}$. %
Definition~\ref{def-w-sol-tp} is silent about the initial data $\textbf{u}_0$ and
the forces field $\textbf{f}$. But, this will be clear when we bellow establish the existence result. %

\subsection{Existence result}\label{sec-existence}

To prove the existence of a weak solution of the modified NS problem (\ref{NSE-2})-(\ref{b-c}), we will adapt the same arguments used to prove existence for the classical NS problem (see, \emph{e.g.} Lions~\cite{Lions-1969}, Temam~\cite{Temam-1979}, Galdi~\cite{Galdi-2000}). However, it is worth to notice that, in (\ref{NSE-2}), additionally to the usual nonlinear term for the classical NS equations, $(\textbf{u$\cdot\nabla$})\textbf{u}$, we have another one, the absorption term
$\alpha|\textbf{u}|^{\sigma-2}\textbf{u}$. %
We shall adapt the proof for the classical NS problem in any dimension $N\geq 2$ given in Lions~\cite[pp. 75-77]{Lions-1969}.

\begin{theorem}\label{existence-ws}
Assume that $\textbf{f}\in\mathrm{L}^2(0,T;\textbf{V}')$ and $\textbf{u}_0\in\textbf{H}$. Then, there exists, at least, a weak solution of the modified NS problem (\ref{NSE-2})-(\ref{b-c}) in the sense of Definition~\ref{def-w-sol-tp}.
\end{theorem}
PROOF. %
\emph{1. Existence of approximate solutions.} %
Let $s\geq 1$ be the smaller integer not lesser than $N/2$. We may assume that $s>1$, because for $s=1$, $\textbf{V}_{s}=\textbf{V}$ and the proof would follow in the same manner, but even simpler. %
We consider the basis $\left\{\textbf{v}_k\right\}_{k\in\textbf{N}}$ of $\textbf{V}_{s}$,  given by the (non-zero) solutions $\textbf{v}_j$ of the following spectral problem associated to the eigenvalues $\lambda_j>0$:
\begin{equation*}\label{spectral-p}
\int_{\Omega}\sum_{|\alpha|=s}\mathrm{D}^{\alpha}\textbf{v}_j\cdot\mathrm{D}^{\alpha}\mathbf{\varphi}\,d\textbf{x}=
\lambda_j
\int_{\Omega}\textbf{v}_j\cdot\mathbf{\varphi}\,d\textbf{x}\quad \forall\ \mathbf{\varphi}\in\textbf{V}_{s}\,.
\end{equation*}
Since $\textbf{V}_s\hookrightarrow\textbf{V}\hookrightarrow\textbf{H}$, $\left\{\textbf{v}_k\right\}_{k\in\textbf{N}}$ can be chosen as being an orthonormal basis in $\textbf{H}$. %
Let us consider also the corresponding $m$-dimensional space, say $\textbf{V}^{m}$, spanned by $\textbf{v}_1$, \dots, $\textbf{v}_m$. %
For each $m\in\mathbb{N}$, we search an approximate solution
$\textbf{u}_m$ of (\ref{w-s-s-h-t}) in the form
\begin{equation}\label{app-sol-t}
\textbf{u}_m=\sum_{k=1}^m c_{km}(t)\textbf{v}_k\,,
\end{equation}
where $\textbf{v}_k\in\textbf{V}^{m}$ and $c_{km}(t)$ are the functions we look for. %
These functions are found by solving the following system of ordinary differential
equations obtained from (\ref{w-s-s-h-t}):
$$\frac{d}{dt}\int_{\Omega}\textbf{u}_m(t)\cdot\textbf{v}_k\,d\textbf{x}+\nu\int_{\Omega}\mathbf{\nabla}(\textbf{u}_m(t)):\mathbf{\nabla}\,\textbf{v}_kd\textbf{x}+
\int_{\Omega}(\textbf{u}_m(t)\mathbf{\cdot\nabla})\textbf{u}_m(t)\cdot\textbf{v}_kd\textbf{x}$$
\begin{equation}\label{app-weak-sol-N2}
\vspace{-0.5cm}
\end{equation}
$$+\alpha\int_{\Omega}|\textbf{u}_m(t)|^{\sigma-2}\textbf{u}_m(t)\cdot\textbf{v}_kd\textbf{x}
=\int_{\Omega}\textbf{f}(t)\cdot\textbf{v}_kd\textbf{x}\,;$$
\begin{equation}\label{app-in-cond}
c_{km}(0)=\int_{\Omega}\textbf{u}_{0m}\cdot\textbf{v}_kd\textbf{x}\,;
\end{equation}
for $k=1,\dots,m$ and where $\textbf{u}_{0m}=\textbf{u}_{m}(0)\in\textbf{V}^m$ is such that
\begin{equation}\label{str-H-u0}
\mbox{$\textbf{u}_{0m}\to\textbf{u}_{0}$ strongly in $\textbf{H}$ as $m\to\infty$.} %
\end{equation}
From the elementary theory of ordinary differential equations, problem (\ref{app-sol-t})-(\ref{app-in-cond}) has a unique solution
$c_{km}\in
\mathrm{C}^1([0,T_m])$, for some small interval of time $[0,T_m]\subset[0,T]$.

\vspace{0.2cm}\noindent
\emph{2. A priori estimates I.} We multiply (\ref{app-weak-sol-N2}) by $c_{km}(t)$ and add these equations from
$k=1$ to $k=m$. %
Then, according to Lions~\cite[Lemme I-6.5]{Lions-1969} (see also Temam~\cite[Lemma II-1.3]{Temam-1979}), we obtain
\begin{equation}\label{apriori-eq}
\begin{split}
\frac{1}{2}\frac{d}{d\,t}\|\textbf{u}_m(t)\|_{\textbf{L}^2(\Omega)}^2&\,+\nu\|\mathbf{\nabla}\textbf{u}_m(t)\|^2_{\textbf{L}^2(\Omega)}
+\alpha\|\textbf{u}_m(t)\|^{\sigma}_{\textbf{L}^{\sigma}(\Omega)}\\
&\,=
\int_{\Omega}\textbf{u}_m(t)\cdot\textbf{f}(t)\,d\textbf{x}\,.
\end{split}
\end{equation}
Integrating (\ref{apriori-eq}) from $0$ to $t\leq T_m$, using Schwarz's and Cauchy's (with a suitable $\varepsilon>0$) inequalities, we achieve to
\begin{equation}\label{apriori-ineq2}
\begin{split}
\|\textbf{u}_m(t)\|_{\textbf{L}^2(\Omega)}^2+\nu\int_0^t\|\mathbf{\nabla}\textbf{u}_m(s)&\|^2_{\textbf{L}^2(\Omega)}ds
+2\alpha\int_0^t\|\textbf{u}_m(s)\|^{\sigma}_{\textbf{L}^{\sigma}(\Omega)}ds
\leq\\
\|\textbf{u}_{0\,m}\|_{\textbf{L}^2(\Omega)}^2+&\frac{1}{\nu}\int_0^t\|\textbf{f}(s)\|_{V'}^2ds
\end{split}
\end{equation}
for $t<T_m$. %
Since $\|\textbf{u}_{0\,m}\|_{\textbf{L}^2(\Omega)}\leq\|\textbf{u}_{0}\|_{\textbf{L}^2(\Omega)}$, the assumptions $\textbf{u}_{0}\in\textbf{H}$ and $\textbf{f}\in\mathrm{L}^2(0,T;\textbf{V}')$ justify that the left-hand side of (\ref{apriori-ineq2}) is finite. In particular, it follows that $|c_{km}(t)|^2<\infty$ for all $k=1,\dots,m$. %
In consequence, by standard results on ordinary differential equations, we get that $T_m=T$ for all $m\in\textbf{N}$,
otherwise $|c_{km}(t)|\to\infty$ as $t\to T_m$. %
Moreover, from (\ref{apriori-ineq2}) and once that $\nu,\ \alpha>0$, we obtain
\begin{equation}\label{apriori-ineq3}
\sup_{t\in[0,T]}\|\textbf{u}_m(t)\|_{\textbf{L}^2(\Omega)}^2\leq
\|\textbf{u}_0\|_{\textbf{L}^2(\Omega)}^2+\frac{1}{\nu}\|\textbf{f}\|_{\textbf{L}^2(0,T;\textbf{V'})}^2\,,
\end{equation}
which implies that
\begin{equation}\label{bound-inf}
\mbox{$\textbf{u}_m$ remains bounded in $\mathrm{L}^{\infty}(0,T;\textbf{H})$}.
\end{equation}
On the other hand, if we replace $t$ by $T$ in
(\ref{apriori-ineq2}), we obtain
\begin{equation}\label{apriori-ineq2-T}
\begin{split}
\|\textbf{u}_m(T)\|_{\textbf{L}^2(\Omega)}^2+\nu&\int_0^T\|\textbf{u}_m(t)\|^2_{\textbf{H}^1(\Omega)}dt+
2\alpha\int_0^T\|\textbf{u}_m(t)\|^{\sigma}_{\textbf{L}^{\sigma}(\Omega)}dt
\leq\\
&\|\textbf{u}_0\|_{\textbf{L}^2(\Omega)}^2+\frac{1}{\nu}\|\textbf{f}\|_{\textbf{L}^2(0,T;\textbf{V'})}^2\,. %
\end{split}
\end{equation}
This estimate enables us to say that
\begin{equation}\label{bound-2-sigma}
\mbox{$\textbf{u}_m$ remains bounded in $\mathrm{L}^{2}(0,T;\textbf{V})$ and in $\mathrm{L}^{\sigma}(Q_T)$.}
\end{equation}
Moreover,
\begin{equation}\label{bound-sigma'}
\mbox{$|\textbf{u}_m|^{\sigma-2}\textbf{u}_m$ remains bounded in $\mathrm{L}^{\sigma'}(Q_T)$.}
\end{equation}

\vspace{0.2cm}\noindent
\emph{3. A priori estimates II.} %
Let us consider the orthogonal projection $P_m:\textbf{H}\to\textbf{V}^m$, $P_m(\textbf{u})=\sum_{k=1}^m\int_{\Omega}\textbf{u}\cdot\textbf{v}_k\,d\textbf{x}\,\textbf{v}_k$.
We thus obtain from (\ref{app-weak-sol-N2})
\begin{equation}\label{Pm-NS}
\frac{\partial\,\textbf{u}_m}{\partial\,t}=P_m\left(\nu\triangle\textbf{u}_m\right)-P_m\left((\textbf{u}_m\cdot\mathbf{\nabla})\textbf{u}_m\right)
-P_m\left(|\textbf{u}_m|^{\sigma-2}\textbf{u}_m\right)+P_m\left(\textbf{f}\right)\,,
\end{equation}
where we have used $P_m(\textbf{u}_m)=\textbf{u}_m$ by virtue of (\ref{app-sol-t}). %
Using (\ref{apriori-ineq3}) and (\ref{apriori-ineq2-T}), and the special choice of the basis of $\textbf{V}^m$, we deduce, arguing as in Lions~\cite[I-6.4.3]{Lions-1969}, that the sequences
$P_m\left(\nu\triangle\textbf{u}_m\right)$ and $P_m\left((\textbf{u}_m\cdot\mathbf{\nabla})\textbf{u}_m\right)$ are bounded in $\mathrm{L}^2(0,T;\textbf{V}'_s)$. %
Moreover, from hypothesis, one clearly has $\textbf{f}\in\mathrm{L}^2(0,T;\textbf{V}'_s)$. %
Then, from (\ref{Pm-NS}),
\begin{equation}\label{bound-ddt}
\mbox{$\displaystyle\frac{\partial\,\textbf{u}_m}{\partial\,t}$ remains bounded in $\mathrm{L}^2(0,T;\textbf{V}'_s)+\textbf{L}^{\sigma'}(Q_T)$.}
\end{equation}
\vspace{0.2cm}\noindent
\emph{4. Passing to the limit.}
From (\ref{bound-inf}), (\ref{bound-2-sigma}), (\ref{bound-sigma'}) and (\ref{bound-ddt}), there exist functions $\mathbf{u}$ and $\mathbf{v}$, and there exists a subsequence, which we still denote by $\textbf{u}_m$, such that
\begin{equation}\label{ws-c-Linf}
\textbf{u}_m\to\textbf{u}\quad\mbox{weak-star in}\quad\mathrm{L}^{\infty}(0,T;\textbf{H})\quad\mbox{as}\quad m\to\infty\,,
\end{equation}
\begin{equation}\label{w-c-L2}
\textbf{u}_m\to\textbf{u}\quad\mbox{weakly in}\quad\mathrm{L}^2(0,T;\textbf{V})\quad\mbox{as}\quad m\to\infty\,,
\end{equation}
\begin{equation}\label{w-c-L-sigma}
\textbf{u}_m\to\textbf{u}\quad\mbox{weakly in}\quad\mathrm{L}^{\sigma}(Q_T)\quad\mbox{as}\quad m\to\infty\,,
\end{equation}
\begin{equation}\label{w-c-L-sigma'}
|\textbf{u}_m|^{\sigma-2}\textbf{u}_m\to\textbf{v}\quad\mbox{weakly in}\quad\mathrm{L}^{\sigma'}(Q_T)\quad\mbox{as}\quad m\to\infty
\end{equation}
and
\begin{equation}\label{w-c-ddt}
\frac{\partial\,\textbf{u}_m}{\partial\,t}\to\frac{\partial\,\textbf{u}}{\partial\,t}\quad\mbox{weakly in}\quad\mathrm{L}^2(0,T;\textbf{V}'_s)+\textbf{L}^{\sigma'}(Q_T)\quad\mbox{as}\quad m\to\infty\,.
\end{equation}
Then, due to (\ref{distr-incl}), (\ref{w-c-L2}) and (\ref{w-c-ddt}), and according to a well known compactness result (see Lions~\cite[Théorème I-5.1]{Lions-1969}),
\begin{equation}\label{s-c-L2}
\textbf{u}_m\to\textbf{u}\quad\mbox{strongly in}\quad\mathrm{L}^{2}(0,T;\textbf{H})\quad\mbox{as}\quad m\to\infty\,.
\end{equation}
Now, we multiply (\ref{app-weak-sol-N2}) by
$\psi\in\mathrm{C}^1([0,T])$, with $\psi(T)=0$, and then we integrate the resulting equations from $0$ to $T$. %
We thus obtain
\begin{equation}\label{app-weak-sol-N2-int-t}
\begin{split}
-\int_0^T\int_{\Omega}\textbf{u}_m(t)\cdot\textbf{v}_k\,\psi'(t)\,d\textbf{x}dt
+\nu&\int_0^T\int_{\Omega}\mathbf{\nabla}\textbf{u}_m(t):\mathbf{\nabla}\,\textbf{v}_k\,\psi(t)\,d\textbf{x}dt+ \\
\int_0^T\int_{\Omega}(\textbf{u}_m(t)\mathbf{\cdot\nabla})\textbf{u}_m(t)\cdot\textbf{v}_k\,\psi(t)\,d\textbf{x}dt&+
\alpha\int_0^T\int_{\Omega}|\textbf{u}_m(t)|^{\sigma-2}\textbf{u}_m(t)\cdot\textbf{v}_k\,\psi(t)d\textbf{x}dt\\
=\int_0^T\int_{\Omega}\textbf{f}(t)\cdot\textbf{v}_k\,\psi(t)\,&d\textbf{x}dt+
\int_{\Omega}\textbf{u}_{m0}\cdot\textbf{v}_k\,\psi(0)\,d\textbf{x}\,.
\end{split}
\end{equation}
Now, we can pass to the limit in (\ref{app-weak-sol-N2-int-t}) by using (\ref{w-c-L2}) in the linear terms (\ref{app-weak-sol-N2-int-t})$_1$ and (\ref{app-weak-sol-N2-int-t})$_2$, (\ref{w-c-L2}) along with (\ref{s-c-L2}) in the nonlinear term (\ref{app-weak-sol-N2-int-t})$_3$ and (\ref{str-H-u0}) in the term (\ref{app-weak-sol-N2-int-t})$_6$ (see Lions~\cite[I-6.4]{Lions-1969}). %
With respect to the term (\ref{app-weak-sol-N2-int-t})$_4$, we notice that by considering a new subsequence, we may assume that $\textbf{u}_n\to\mathbf{u}$ a.e. in $Q_T$. %
This and (\ref{w-c-L-sigma'}) prove that $\textbf{v}=|\textbf{u}|^{\sigma-2}\textbf{u}$, and we can pass to the limit in this term as well. %
Finally, using linear and continuity arguments, we can prove that $\textbf{u}$ satisfies
(\ref{w-s-s-h-t}) in the distribution sense of $\textrm{C}_0^{\infty}(0,T)$. %
Moreover, using standard arguments, we can prove that $\textbf{u}$ satisfies 2. of Definition~\ref{def-w-sol-tp}. See \emph{e.g.} Temam~\cite[pp. 288-289]{Temam-1979} for the details.$\square$

\begin{remarks}
1. This existence result still holds for unbounded domains $\Omega$ as far as the used Sobolev type inequalities hold. The main difference lies in the fact that the imbedding $\textbf{V}\hookrightarrow\textbf{H}$ is no longer compact. But this difficulty is overcame in a standard manner (see \emph{e.g.} Temam~\cite[Remark III-3.2]{Temam-1979}).\\
2. From what we have said at the beginning of this section, in particular (\ref{distr-incl}), the assumption $\textbf{f}\in\mathrm{L}^2(0,T;\textbf{V}')$
can be replaced by $\textbf{f}\in\mathrm{L}^2(0,T;\textbf{H}^{-1}(\Omega))$\,.
On the other hand, with some minor modifications in the proof, we can assume that $\textbf{f}\in\mathrm{L}^2(0,T;\textbf{L}^p(\Omega))$ for some some $p\geq 1$.
\end{remarks}

\subsection{Energy relation}\label{sec-en-ineq}

In this section we shall establish the energy relations satisfied by the weak solutions of the modified NS problem
(\ref{NSE-2})-(\ref{b-c}). %
Before we establish the main result, let us define the (kinetic) energy associated with the problem (\ref{NSE-2})-(\ref{b-c}):
\begin{equation}\label{k-energy}
E(t):=\frac{1}{2}\|\textbf{u}(t)\|_{\textbf{L}^2(\Omega)}^2\,.
\end{equation}
As one should expect, at least from the physical point of view, any solution of problem (\ref{NSE-2})-(\ref{b-c}) should be such that the associated kinetic energy at a certain time $t$ (considering $\textbf{f}\equiv\textbf{0}$) is equal to $E(s)$, for some $0\leq s\leq t$, minus the energy dissipated by viscosity in $[s,t]$,
$$\nu\int_s^t\|\mathbf{\nabla\,u}(\tau)\|^2_{\textbf{L}^2(\Omega)}d\tau\quad \mbox{(enstrophy)}\,,$$
and minus the energy dissipated by a sink associated to $\alpha$ in $[s,t]$,
$$\alpha\int_s^t\|\textbf{u}(\tau)\|^{\sigma}_{\textbf{L}^{\sigma}(\Omega)}d\tau\quad \mbox{(absorption)}\,.$$
However, even for the classical NS problem, weak solutions obey to such an equality only if $N=2$. For $N\geq 3$, it is only possible to prove they satisfy to an energy inequality (see \emph{e.g.} Galdi~\cite{Galdi-2000}). %
On the other hand, from the above considerations, we should expect
\begin{equation}\label{u0-en-fin}
\textbf{u}_0\in\textbf{H} \Rightarrow E(0)<\infty \Rightarrow E(t)\leq E(0)<\infty\ \mbox{for all}\ t\geq 0\,.
\end{equation}
In the following result we shall see the weak solutions of the modified NS problem (\ref{NSE-2})-(\ref{b-c}) satisfy to an energy inequality, regardless the domain dimension.

\begin{theorem}\label{Th-energy-in}
Assume that $\textbf{f}\in\mathrm{L}^2(0,T;\textbf{V}')$ and $\textbf{u}_0\in\textbf{H}$.
Let $\textbf{u}$ be a weak solution of the problem (\ref{NSE-2})-(\ref{b-c}) in the sense of Definition~\ref{def-w-sol-tp}.
Then this solution satisfies to
\begin{equation}\label{energy-in-df}
\begin{split}
\frac{1}{2}\frac{d}{d\,t}\|\textbf{u}(t)\|^2_{\textbf{L}^2(\Omega)}+ & \nu\|\mathbf{\nabla\,u}(t)\|^2_{\textbf{L}^2(\Omega)}+
 \alpha\|\textbf{u}(t)\|^{\sigma}_{\textbf{L}^{\sigma}(\Omega)}\leq \\
& \int_{\Omega}\textbf{u}(t)\cdot\textbf{f}(t)\,d\textbf{x}
\end{split}
\end{equation}
for a.a. $t\in[0,T]$, and
\begin{equation}\label{energy-in}
\begin{split}
\|\textbf{u}(t)\|^2_{\textbf{L}^2(\Omega)}+ & 2\nu\int_0^t\|\mathbf{\nabla\,u}(s)\|^2_{\textbf{L}^2(\Omega)}ds+
2\alpha\int_0^t\|\textbf{u}(s)\|^{\sigma}_{\textbf{L}^{\sigma}(\Omega)}ds\leq \\
& 2\int_0^t\int_{\Omega}\textbf{u}(s)\cdot\textbf{f}(s)\,d\textbf{x}ds+\|\textbf{u}_0\|^2_{\textbf{L}^2(\Omega)}
\end{split}
\end{equation}
for a.a. $t\in[0,T]$.
\end{theorem}
PROOF. To prove (\ref{energy-in-df}) we take the limit\,~inf, as $m\to\infty$, of (\ref{apriori-eq}). Then from
(\ref{str-H-u0}), (\ref{ws-c-Linf})-(\ref{w-c-L-sigma}) and a classical property of weak limits (see Galdi~\cite{Galdi-2000}), the result follows.
For (\ref{energy-in}), we start by integrating (\ref{apriori-eq}) between $0$ and $t\leq T$. We thus obtain
\begin{equation}\label{energy-in-apriori}
\begin{split}
\|\textbf{u}_m(t)\|_{\textbf{L}^2(\Omega)}^2+ & 2\nu\int_0^t\|\mathbf{\nabla}(\textbf{u}_m(s))\|^2_{\textbf{L}^2(\Omega)}ds
+ 2\alpha\int_0^t\|\textbf{u}_m(s)\|^{\sigma}_{\textbf{L}^{\sigma}(\Omega)}ds=\\
& \int_0^t\int_{\Omega}\textbf{u}_m(s)\cdot\textbf{f}(s)\,ds+\|\textbf{u}_{m0}\|_{\textbf{L}^2(\Omega)}^2\,.
\end{split}
\end{equation}
Now, we take the limit\,~inf, as $m\to\infty$, of (\ref{energy-in-apriori}). Then from
(\ref{str-H-u0}), (\ref{ws-c-Linf})-(\ref{w-c-L-sigma}) and the same classical property of weak limits, it follows (\ref{energy-in}). $\square$

\bigskip\noindent
Such as for the classical NS problem, we can prove the weak solutions of the modified NS problem
(\ref{NSE-2})-(\ref{b-c}) belonging to $\mathrm{L}^4(0,T;\textbf{L}^4(\Omega))$ satisfy to the energy equality
$$\|\textbf{u}(t)\|^2_{\textbf{L}^2(\Omega)}+2\nu\int_0^t\|\mathbf{\nabla\,u}(s)\|^2_{\textbf{L}^2(\Omega)}ds+
2\alpha\int_0^t\|\textbf{u}(s)\|^{\sigma}_{\textbf{L}^{\sigma}(\Omega)}ds=$$
\begin{equation}\label{energy-eq}
\vspace{-0.5cm}
\end{equation}
$$2\int_0^t\int_{\Omega}\textbf{u}(s)\cdot\textbf{f}(s)\,ds+\|\textbf{u}_0\|^2_{\textbf{L}^2(\Omega)}$$
for a.a. $t\in[0,T]$. The proof can be adapted from Galdi~\cite[Theorem 4.1]{Galdi-2000}. %
The only difference lies in the absorption term which can be handle such as in the proof of Theorem~\ref{existence-ws}. Moreover, if $N=2$, every weak solution of the modified NS problem satisfies to the energy equality (\ref{energy-eq}).
Indeed for $N=2$ any such weak solution $\textbf{u}$ satisfies to the so-called Ladyzhenskaya inequality
$$\int_0^T\|\textbf{u}(t)\|_{\textbf{L}^4(\Omega)}dt\leq C
  \int_0^T\|\textbf{u}(t)\|_{\textbf{L}^2(\Omega)}^2\|\mathbf{\nabla\,u}(t)\|_{\textbf{L}^2(\Omega)}^2dt<\infty\,.$$
This is no longer valid for $N\geq 3$ and therefore, such as for the classical Navier-Stokes problem, the question of whether a weak solution of the modified NS problem (\ref{NSE-2})-(\ref{b-c}) obeys the energy equality  (\ref{energy-eq}) remains open. %

\subsection{Uniqueness result}\label{sec-uniqueness}

In the mathematical theory of the NS equations another important open problem is the uniqueness of weak solutions for $N\geq 3$. It is possible to show uniqueness for all time under additional conditions on the admissible weak solutions, or in a small interval of time by assuming small data compared to viscosity (see \emph{e.g.} Galdi~\cite{Galdi-2000}). However, for $N=2$ the mathematical theory is complete and the existence result is now classical.
In consequence, we will show the 2-D modified NS problem (\ref{NSE-2})-(\ref{b-c}) inherits this property. The proof is again an adaptation of the classical result (see \emph{e.g.} Lions~\cite{Lions-1969}) and the crucial part lies in the absorption term.

\begin{theorem}\label{theorem-unique-N2}
Assume that $\textbf{f}\in\mathrm{L}^2(0,T;\textbf{V}')$ and $\textbf{u}_0\in\textbf{H}$.
Let $\textbf{v}$ and $\textbf{w}$ be two weak solutions of the modified NS problem (\ref{NSE-2})-(\ref{b-c}) in the sense of Definition~\ref{def-w-sol-tp}. Then $\textbf{w}=\textbf{v}$ a.e. in $Q_T$.
\end{theorem}
PROOF. Arguing as in Temam~\cite[Theorem
III-3.2]{Temam-1979}, we get from (\ref{w-s-s-h-t}) the following relation for
$\textbf{u}=\textbf{v}-\textbf{w}$:
\begin{equation}\label{energy-1}
\frac{d}{dt}\|\textbf{u}(t)\|_{\textbf{L}^2(\Omega)}^2+2\nu\|\mathbf{\nabla}\textbf{u}(t)\|_{\textbf{L}^2(\Omega)}^2+I_1=
I_2\,,
\end{equation}
where
$$I_1:=2\alpha\int_{\Omega}\left(|\textbf{v}(t)|^{\sigma-2}\textbf{v}(t)-|\textbf{w}(t)|^{\sigma-2}\textbf{w}(t)\right)\cdot\textbf{u}(t)\,d\textbf{x}\,,$$
$$I_2:=2\int_{\Omega}\left[(\textbf{w}(t)\mathbf{\cdot\nabla})\textbf{w}(t)-(\textbf{v}(t)\mathbf{\cdot\nabla})\textbf{v}(t)\right]\cdot\textbf{u}(t)\,d\textbf{x}\,.$$%
In the term $I_1$, we use (\ref{gin-xi-eta}) with $\mathbf{\xi}=\textbf{v}$, $\mathbf{\eta}=\textbf{w}$, $p=\sigma$ and $\delta=0$, to prove that $I_1\geq 0$. %
On the other hand, it can be proved, in a standard manner, that
\begin{equation}\label{est-I2}
|I_2|\leq
\frac{2\nu}{C}\|\textbf{u}(t)\|_{\textbf{H}^1(\Omega)}^2
+\frac{C}{\nu}\|\textbf{u}(t)\|_{\textbf{L}^2(\Omega)}^2\|\textbf{w}(t)\|_{\textbf{H}^1(\Omega)}^2\,,
\end{equation}
for some positive constant. %
Then, using Poincaré's inequality, we obtain, from
(\ref{energy-1}) and (\ref{est-I2}), the following relation
\begin{equation}\label{energy-2}
\frac{d}{dt}\|\textbf{u}(t)\|_{\textbf{L}^2(\Omega)}^2\leq
\frac{C}{\nu}\|\textbf{u}(t)\|_{\textbf{L}^2(\Omega)}^2\|\textbf{w}(t)\|_{\textbf{H}^1(\Omega)}^2\,.
\end{equation}
Integrating (\ref{energy-2}), using (\ref{apriori-ineq2-T}) for
$\textbf{w}$, and known that $\textbf{u}(0)=\textbf{0}$, we prove
that $\textbf{v}=\textbf{w}$ a.e. in $Q_T$.$\square$

\bigskip\noindent
For $N\geq 3$ and analogously as to the classical NS problem, one can proves the weak solution of the modified NS problem, satisfying (\ref{energy-in}), is unique in the class $\mathrm{L}^2(0,T;\textbf{V})\cap\mathrm{L}^{\infty}(0,T;\textbf{H})\cap\mathrm{L}^{\sigma}(Q_T)\cap\mathrm{L}^{s}(0,T;\textbf{L}^r(\Omega))$ with $2/s+N/r\leq 1$ and for $r>N$ (see Lions~\cite[Théorème I-6.9]{Lions-1969}). In particular the uniqueness holds in the class
$\mathrm{L}^2(0,T;\textbf{V})\cap\mathrm{L}^{\infty}(0,T;\textbf{H})\cap\mathrm{L}^{\sigma}(Q_T)$ for $\sigma\geq N+2$.

\section{Asymptotic stability}

In this section we shall study the behavior in time of the weak solutions of the modified NS problem (\ref{NSE-2})-(\ref{b-c}). %
Here, we shall assume that $T$ is sufficiently large or even let $T=\infty$. %
It is worth to recall that, to the best of our knowledge, the late studies on the asymptotic behavior of the weak solutions of the classical NS problem provide only power-time decays (see the references cited in Section~\ref{sect-TMNS}). %
The first result we present here, establishes a temporal qualitative property satisfied by the weak solutions of the modified NS problem (\ref{NSE-2})-(\ref{b-c}), which is usually referred to as the extinction in (a finite) time property. %

\begin{theorem}[Extinction in time]\label{main-theorem}
Assume $1<\sigma<2$ and $\textbf{u}_0\in\textbf{H}$, and let $\textbf{u}$ be a weak solution of the modified NS problem
(\ref{NSE-2})-(\ref{b-c}) in the sense of Definition~\ref{def-w-sol-tp}. %
\begin{enumerate}
      \item If $\textbf{f}=\textbf{0}$ a.e. in $Q_T$,  then there exists
$t^{\ast}>0$ such that $\textbf{u}=\textbf{0}$ a.e. in $\Omega$ and for almost all $t\geq t^\ast$.
      \item If $\textbf{f}\not=\textbf{0}$, but
\begin{equation}\label{cond-f-+}
\|\textbf{f}(t)\|_{\textbf{V}'} \leq \epsilon %
\left(1-\frac{t}{{t}_{\textbf{f}}}\right)_{+}^{\frac{1}{2(\mu-1)}}\quad
\mbox{for a.a. $t\in[0,T]$}\,,
\end{equation}
where $t_{\textbf{f}}$ is a fixed positive time and $\mu$ is given by (\ref{mu}),
then there exists a constant $\epsilon_{0}>0$ such
that $\textbf{u}=\textbf{0}$ a.e. in $\Omega$ and for almost all $t\geq t_{\textbf{f}}$, provided $0<\epsilon\leq\epsilon_{0}$.
\end{enumerate}
\end{theorem}
PROOF. \emph{First Step.} If $\textbf{f}=\textbf{0}$ a.e. in $Q_T$, we obtain from (\ref{energy-in-df})
\begin{equation} \label{ei:only11}
 \frac{d}{dt}E(t)+CE_{2,\sigma}(t)\leq 0\quad\mbox{for a.a.}\ t\in[0,T]\,,
\end{equation}
where $E(t)$ is the (kinetic) energy defined in (\ref{k-energy}), $C=\min(\nu,\alpha)$ and
\begin{equation}\label{E-2-sigma}
E_{2,\sigma}(t):=\|\mathbf{\nabla}\textbf{u}(t)\|^{2}_{\textbf{L}^2(\Omega)}+\|\textbf{u}(t)\|^{\sigma}_{\textbf{L}^{\sigma}(\Omega)}\,.
\end{equation}
Using a vectorial version of the Gagliardo-Nirenberg-Sobolev inequality (\ref{Nirenberg}) with $p=q=2$ and  $r=\sigma$, \begin{equation}\label{eq:aux-2p}%
\|\textbf{u}(t)\|_{\textbf{L}^2(\Omega)}\leq C_{GNS}
\|\mathbf{\nabla}\textbf{u}(t)\|_{\textbf{L}^2(\Omega)}^{\theta}
\|\textbf{u}(t)\|_{\textbf{L}^{\sigma}(\Omega)}^{1-\theta}
\quad\mbox{for a.a.}\ t\in[0,T]\,,
\end{equation}%
where $C_{GNS}=C(N,\sigma)$ is the constant resulting from applying (\ref{Nirenberg}) and, according to (\ref{theta}),
\begin{equation}\label{theta-s}
\theta=1-\frac{2\sigma}{(2-\sigma)N+2\sigma}\,.
\end{equation}
Then, using the algebraic inequality (\ref{algebraic}), we obtain
\begin{equation}\label{eq:aux-2}%
\|\textbf{u}(t)\|_{\textbf{L}^2(\Omega)}^{2}\leq C_{GNS}^2 E_{2,\sigma}(t)^{\mu}\quad\mbox{for a.a.}\ t\in[0,T]\,,
\end{equation}%
where, from (\ref{eq:aux-2p}) and (\ref{theta-s}),
\begin{equation}\label{mu}
\mu=1+\frac{2(2-\sigma)}{(2-\sigma)N+2\sigma}\,.
\end{equation}
Then (\ref{ei:only11}) and (\ref{eq:aux-2}) lead us to the homogeneous ordinary differential inequality
\begin{equation} \label{odi}
\frac{d}{dt}E(t)+CE(t)^{1/\mu}\leq 0\quad\mbox{for a.a.}\ t\in[0,T]\,,
\end{equation}
where now $C=\min(\nu,\alpha)(2/C_{GNS}^2)^{1/\mu}$. %

\vspace{0.15cm}\noindent \emph{Second Step.} %
In order to integrate (\ref{odi}), we need to analyze the exponent of nonlinearity $\mu$ given by (\ref{mu}). %
But, firstly, we recall that $\mu$ is written in terms of the interpolation exponent $\theta$. %
According to Lemma~\ref{L-Nirenberg} with $p=q=2$ and $r=\sigma$, and also (\ref{theta-s}), the admissible range of $\theta$ shows us that
$0\leq\theta\leq 1$ iff $0\leq \sigma \leq 2$ for $N=1$ or $N\geq 3$,
and $0\leq\theta<\frac{2}{2+\sigma}$ iff $0< \sigma \leq 2$ for $N=2$. %
In consequence, within these values of $\sigma$, we see that $\mu>1$ iff $1<\sigma<2$. %
The other possible values of $\sigma$, \emph{i.e.} $0\leq\sigma\leq 1$ go out of our initial assumption that $\sigma>1$.

\vspace{0.15cm}\noindent \emph{Third Step.} %
An explicit integration of (\ref{odi}) between $t=0$ and $t$ lead us to
\begin{equation}\label{in-E-E0}
E(t)\leq\left(E(0)^{\frac{\mu-1}{\mu}}-Ct\right)^{\frac{\mu}{\mu-1}}\,,\quad C=\frac{\mu-1}{\mu}\min(\nu,\alpha)\left(\frac{2}{C_{GNS}^{2}}\right)^{\frac{1}{\mu}}\,.
\end{equation}
Notice that, from (\ref{u0-en-fin}), the initial energy $E(0)$ is finite, and
\begin{equation}\label{frac-mu}
\mu>1\Rightarrow 0<\frac{\mu}{\mu-1}=\frac{4+(2-\sigma)N}{2(2-\sigma)}\,.
\end{equation}
Using (\ref{k-energy}), (\ref{mu}) and (\ref{frac-mu}), we prove that the right-hand side of (\ref{in-E-E0}) vanishes for
\begin{equation}\label{t-ast}
t=\frac{E(0)^{\frac{\mu-1}{\mu}}}{C}
\equiv\|u_0\|_{\textbf{L}^2(\Omega)}^{\frac{4(2-\sigma)}{4+(2-\sigma)N}}
\frac{4+(2-\sigma)N}{4(2-\sigma)}\frac{C_{GNS}^{\frac{2\left[(2-\sigma)N+2\sigma\right]}{4+(2-\sigma)N}}}{\min(\nu,\alpha)}\,,
\end{equation}
which proves the first assertion with $t^{\ast}$ given by (\ref{t-ast}). %

\vspace{0.15cm}\noindent \emph{Fourth Step.}
If $\textbf{f}\not=\textbf{0}$, we use Schwarz's and Cauchy's inequalities, the later with $\varepsilon=\nu/2$, to obtain for a.a. $t\in[0,T]$
\begin{equation}\label{ineq-121}
\left|\int_{\Omega}\textbf{u}(t)\cdot\textbf{f}(t)\,d\textbf{x}\right|\leq
\|\textbf{u}(t)\|_{V}\,||\textbf{f}(t)||_{V'}
\leq \frac{\nu}{2} \|\mathbf{\nabla}\textbf{u}(t)\|_{\textbf{L}^2(\Omega)}^2+
\frac{1}{2\nu}\|\textbf{f}(t)\|_{\textbf{V}'}^2.
\end{equation}
Then, from (\ref{energy-in-df}) and (\ref{ineq-121}),
\begin{equation} \label{odi-c}
\frac{d}{dt}E(t)+C_1E_{2,\sigma}(t)\leq C_2\|\textbf{f}(t)\|_{\textbf{V}'}^2\quad\mbox{for a.a.}\ t\in[0,T]\,,
\end{equation}
where $C_1=\min(\nu/2,\alpha)$ and $C_2=1/(2\nu)$. %
Using (\ref{cond-f-+}) and (\ref{eq:aux-2}), we obtain the following non-homogeneous ordinary differential inequality
\begin{equation}\label{ineq-15}
\frac{d}{dt}E(t)+C_1E(t)^{1/\mu}\leq C_2\epsilon^2\left(1-\frac{t}{{t}_{\textbf{f}}}\right)_{+}^{\frac{1}{\mu-1}}\qquad\mbox{for a.a. $t\in [0,T]$}\,,
\end{equation}
where $C_1=\frac{\min(\nu/2,\alpha)}{(C_{GNS}^2/2)^{1/\mu}}$ and $C_2=1/(2\nu)$. %
To analyze (\ref{ineq-15}), we need the following result (see the proof in Antontsev \emph{et al.} \cite[\S 1.2.2]{ADS-2002}).
\begin{lemma}
Let $\delta>0$ such that $\left(t_{\textbf{f}}-\delta,t_{\textbf{f}}+\delta\right)\subset[0,T]$ and assume $E\in\mathrm{W}^{1,1}\left(t_{\textbf{f}}-\delta,t_{\textbf{f}}+\delta\right)$ satisfies the differential inequality
\begin{equation}\label{diff-ineq}
\frac{d}{dt}E(t)+\varphi(E(t))\leq F\left(\left(1-\frac{t}{{t}_{\textbf{f}}}\right)_{+}\right)\quad\mbox{a.e. in  $\left(t_{\textbf{f}}-\delta,t_{\textbf{f}}+\delta\right)$}\,,
\end{equation}
where $\varphi$ is a continuous non-decreasing function such that
\begin{equation}\label{cond-varphi}
\varphi(0)=0\quad\mbox{and}\quad \int_{0^+}\frac{ds}{\varphi(s)}<\infty\,,
\end{equation}
and the function $F$ satisfies, for some $\overline{k}\in(0,1)$, to
\begin{equation}\label{F(s)}
F(s)\leq (1-\overline{k})\varphi(\eta_{\overline{k}}(s))\quad\mbox{in $(0,t_{\textbf{f}})$}\,,
\end{equation}
where
$$\eta_{k}(s)=\theta^{-1}_k(s)\quad\mbox{and}\quad \theta_k(s)=\int_0^s\frac{d\tau}{k\varphi(\tau)}\,.$$
Then $E(t)=0$ for all $t\geq t_{\textbf{f}}$.
\end{lemma}
For the relations (\ref{ineq-15}) and (\ref{diff-ineq}) to coincide, let us set
$$\varphi(s)=C_1s^{\frac{1}{\mu}}\quad\mbox{and}\quad F(s)=C_2\epsilon^2s^{\frac{1}{\mu-1}}\,.$$
Then clearly (\ref{cond-varphi}) is satisfied and we have
$$\theta_k(s)=\frac{\mu}{k(\mu-1)}s^{\frac{\mu-1}{\mu}}\quad\mbox{and}\quad
  \eta_k(s)=\left(\frac{k(\mu-1)}{\mu}s\right)^{\frac{\mu}{\mu-1}}\,.$$
Moreover, (\ref{F(s)}) is satisfied if
\begin{equation}\label{epsil-0}
\begin{split}
\epsilon \leq & \sqrt{\frac{C_1}{C_2}(1-\overline{k})\left(\frac{\mu-1}{\mu}\overline{k}\right)^{\frac{1}{\mu-1}}}\\
          = & \sqrt{\frac{2\nu\min(\nu/2,\alpha)}{(C_{GNS}^2/2)^{\frac{(2-\sigma)N+2\sigma}{4+(2-\sigma)N}}}(1-\overline{k})
\left[\frac{2(2-\sigma)}{4+(2-\sigma)N}\overline{k}\right]^{\frac{(2-\sigma)N+2\sigma}{2(2-\sigma)}}}\,.
\end{split}
\end{equation}
Second assertion is thus proved and $\epsilon_0$ is the term defined by the right-hand side of (\ref{epsil-0}).$\square$

\begin{remarks}\label{rem-main-sec}
1. The results established in Theorem~\ref{main-theorem} still hold for unbounded domains $\Omega$ as far as the used Sobolev type inequalities hold (see Lemma~\ref{L-Nirenberg}). %
It is important that the domain is convex and bounded, at least, in one direction.\\
2. We could also have considered non-homogeneous boundary conditions,
say $\mathbf{u}_{B}$, on $\Gamma_T$. But then, in order to carry out
the results of Theorem~\ref{main-theorem}, we would have to assume the
existence of a $t_{B}>0$ such that
$\mathbf{u}_{B}=\mathbf{0}$ for all $t\geq t_{B}$ and
$E(t_{B})<\infty$. In the above proof we only would have to
replace the time $t=0$ by $t=t_{B}$ and, for the second assertion, to presume that $t_{\textbf{f}}>t_{B}$.
\end{remarks}

If we assume, in Theorem~\ref{main-theorem}, $\textbf{f}\in\mathrm{L}^2(0,T;\textbf{L}^p(\Omega))$, then we must replace (\ref{cond-f-+}) by
\begin{equation*}
\|\textbf{f}(t)\|_{\textbf{L}^p(\Omega)} \leq \epsilon %
\left(1-\frac{t}{{t}_{\textbf{f}}}\right)_{+}^{\frac{1}{2(\mu-1)}},\qquad
p=\frac{2N}{N+2},\quad N\not=2\,,
\end{equation*}
where now $\mu$ depends on $p$, $N$ and $\sigma$. Indeed, proceeding such as for (\ref{ineq-121}) and using, in addition, (\ref{Nirenberg}) with $q=p'$, $p=2$ and $\theta=1$, we obtain
\begin{equation*}
\begin{split}
 \left|\int_{\Omega}\textbf{u}(t)\cdot\textbf{f}(t)\,d\textbf{x}\right|\leq & %
  \|\textbf{u}(t)\|_{\textbf{L}^{p'}(\Omega)}\,||\textbf{f}(t)||_{\textbf{L}^{p}(\Omega)}\leq \\
 C\|\mathbf{\nabla\,u}(t)\|_{\textbf{L}^{2}(\Omega)}\,||\textbf{f}(t)||_{\textbf{L}^{p}(\Omega)}
& \leq \varepsilon \|\mathbf{\nabla}\textbf{u}(t)\|_{\textbf{L}^2(\Omega)}^2+
C(\varepsilon)\|\textbf{f}(t)\|_{\textbf{L}^{p}(\Omega)}^2\,.
\end{split}
\end{equation*}
For $N=2$, we assume $\textbf{f}\in\mathrm{L}^2(0,T;\textbf{L}^2(\Omega))$ and by using Poincaré's inequality (\ref{Nirenberg}), we get
\begin{equation*}\label{ineq-12}
\begin{split}
 \left|\int_{\Omega}\textbf{u}(t)\cdot\textbf{f}(t)\,d\textbf{x}\right|\leq & %
  \|\textbf{u}(t)\|_{\textbf{L}^{2}(\Omega)}\,||\textbf{f}(t)||_{\textbf{L}^{2}(\Omega)}\leq \\
 C\|\mathbf{\nabla\,u}(t)\|_{\textbf{L}^{2}(\Omega)}\,||\textbf{f}(t)||_{\textbf{L}^{2}(\Omega)}
& \leq \varepsilon \|\mathbf{\nabla}\textbf{u}(t)\|_{\textbf{L}^2(\Omega)}^2+
C(\varepsilon)\|\textbf{f}(t)\|_{\textbf{L}^{2}(\Omega)}^2\,.
\end{split}
\end{equation*}

The results of Theorem~\ref{main-theorem} can be extended to the limit case of $\sigma=1$. %
In fact, if $\sigma=1$ in (\ref{mu}), then $\mu=1+2/(N+2)$ and all the reasoning of the above proof can be applied.
On the other hand, if $\sigma=2$ these results are no longer valid, because, from (\ref{t-ast}) and (\ref{epsil-0}),
$t^{\ast}\to\infty$ and $\epsilon_0\to 0$ as $\sigma\to 2$. %
In this case, taking $\sigma=2$ in (\ref{mu}), we obtain from (\ref{odi})
\begin{equation}\label{odi-ed}
\frac{d}{dt}E(t)+CE(t)\leq 0\quad\mbox{for a.a.}\ t\in[0,T]\,,\quad C=2\min(\nu,\alpha)/C_{GNS}^2\,.
\end{equation}
Notice that $E_{2,\sigma}(t)\leq E(t)$ is trivial for $\sigma=2$. %
Integrating (\ref{odi-ed}) between $t=0$ and $t>0$, we prove the weak solutions of the modified NS problem (\ref{NSE-2})-(\ref{b-c}), with $\sigma=2$ and $\textbf{f}=\mathbf{0}$, have the following exponential decay
\begin{equation}\label{e-decay-s=2}
\|\textbf{u}(t)\|_{\textbf{L}^2(\Omega)}\leq C_1e^{-C_2t}\quad\mbox{for a.a. $t\geq 0$}\,,
\end{equation}
where $C_1=\|\textbf{u}_0\|_{\textbf{L}^2(\Omega)}$ and
$C_2=2\min(\nu,\alpha)/C_{GNS}^2$. %
The following theorem shows us that for $\sigma>2$, the weak solutions of the modified NS problem (\ref{NSE-2})-(\ref{b-c}), with $\textbf{f}=\mathbf{0}$, have a power-time decay.

\begin{theorem}[Power decay] \label{main-theorem-od}
Assume $\textbf{f}=\mathbf{0}$ and $\textbf{u}_0\in\textbf{H}$, and let $\textbf{u}$ be a weak solution of the modified NS problem
(\ref{NSE-2})-(\ref{b-c}) in the sense of Definition~\ref{def-w-sol-tp}.
If $\sigma>2$, then there exist positive constants $C_1$ and $C_2$ such that
$$\|\textbf{u}(t)\|_{\textbf{L}^2(\Omega)}\leq \left(C_1t+C_2\right)^{-\frac{\sigma+2}{\sigma-2}}\quad\mbox{for a.a. $t\geq 0$}\,.$$
\end{theorem}
PROOF.
Firstly we observe that using Hölder's, Sobolev's (\ref{Nirenberg}) and algebraic (\ref{algebraic}) inequalities, we obatin
\begin{equation}\label{hsa}
\|\textbf{u}(t)\|_{\textbf{L}^2(\Omega)}^2\leq C_S\left(\int_{\Omega}\left(|\textbf{u}(t)|^{\sigma}+|\mathbf{\nabla}\textbf{u}(t)|^{2}\right)\,d\mathbf{x}\right)^{\frac{\sigma+2}{2\sigma}}
\quad\mbox{for a.a.}\ t\in[0,T]\,,
\end{equation}
where $C_S$ is the constant resulting from applying Sobolev's inequality (\ref{Nirenberg}). %
Then (\ref{ei:only11}) and (\ref{hsa}) lead us to the homogeneous ordinary differential inequality
\begin{equation} \label{odi-s>2}
\frac{d}{dt}E(t)+CE(t)^{\frac{2\sigma}{\sigma+2}}\leq 0\quad\mbox{for a.a.}\ t\in[0,T]\,,
\end{equation}
where $C=(2/C_s)^{2\sigma/(\sigma+2)}\min(\nu,\alpha)$. %
Integrating (\ref{odi-s>2}) between $t=0$ and $t>0$ ends the proof, where
$$ C_1=\left(\frac{2}{C_S}\right)^{\frac{2\sigma}{\sigma+2}}\min(\nu,\alpha)\frac{\sigma+2}{\sigma-2}\quad
\mbox{and}\quad
C_2=\left(\frac{\|\textbf{u}_0\|_{\textbf{L}^2(\Omega)}^2}{2}\right)^{-\frac{\sigma-2}{\sigma+2}}\,.\square$$

\begin{remarks}
1. By virtue of using Sobolev's inequality, (\ref{hsa}) holds for any $\sigma\geq\frac{2N}{N+2}$ if $N\geq 3$ or for any $\sigma>1$ if $N=2$, and the assumption $\sigma>2$ satisfies these requirements. \\ %
2. In the limit case $\sigma=2$, (\ref{odi-s>2}) becomes a linear differential inequality analogous to (\ref{odi-ed}) and, again, we obtain the exponential decay (\ref{e-decay-s=2}), where we only have to replace $C_{GNS}^2$ by $C_{S}$.

\end{remarks}

In Theorem~\ref{main-theorem} we have seen that for $1<\sigma<2$ and $\textbf{f}=\textbf{0}$ or for $\textbf{f}$ satisfying to (\ref{cond-f-+}) it was possible to establish that the weak solutions of the modified NS problem
(\ref{NSE-2})-(\ref{b-c}) extinct in a finite time. %
If $\sigma\geq 2$ and $\textbf{f}=\textbf{0}$ it was possible to prove, in Theorem~\ref{main-theorem-od}, that the weak solutions exponentially decay if $\sigma=2$ and decay at a power rate if $\sigma>2$. %
Therefore it is reasonable to ask what happens if $\textbf{f}\not=\textbf{0}$ merely belongs to a suitable function space. %
In order to work this case, let us go back to (\ref{energy-in-df}). Using there Cauchy's inequality with $\epsilon=\nu/2$, we obtain
\begin{equation} \label{ei:only-C-f}
 \frac{d}{dt}E(t)+C_1E_{2,\sigma}(t)\leq C_2\int_{\Omega}|\textbf{f}(t)|^2\,d\textbf{x}\quad\mbox{for a.a.}\ t\in[0,T]\,,
\end{equation}
where $E_{2,\sigma}(t)$ is defined in (\ref{E-2-sigma}) and now $C=\min(\nu/2,\alpha)$ and $C_2=1/(2\nu)$. %
Now we assume that
\begin{equation}\label{C-f}
\int_{\Omega}|\textbf{f}(t)|^2\,d\textbf{x}\leq C_{\textbf{f}}\quad\mbox{for a.a. $t\in[0,T]$}\,,
\end{equation}
where $C_{\textbf{f}}$ is a positive constant. %
If $1<\sigma<2$ we use (\ref{eq:aux-2}) and if $\sigma\geq 2$ we use (\ref{hsa}). %
In any case, we obtain from (\ref{ei:only-C-f}) and (\ref{C-f})
\begin{equation}\label{n.n.o.d.i}
\frac{d}{d\,t}E(t)+C_1E(t)^{\gamma}\leq C_2\quad\mbox{for a.a.}\ t\in[0,T]\,,
\end{equation}
where
\begin{equation}\label{gamma-sigma}
\gamma=\frac{(2-\sigma)N+2\sigma}{4+(2-\sigma)N}\quad\mbox{if}\quad 1<\sigma<2\quad\mbox{and}\quad
\gamma=\frac{2\sigma}{\sigma+2}\quad\mbox{if}\quad \sigma> 2\,,
\end{equation}
\begin{equation}\label{C1C2-sigma}
C_1=\left(\frac{2}{C}\right)^{\gamma}\min\left(\frac{\nu}{2},\alpha\right)\quad\mbox{and}\quad
C_2=\frac{C_{\textbf{f}}}{2\nu}\quad\mbox{for any $\sigma>1$}\,,
\end{equation}
and where $C=C^2_{GNS}$ or $C=C_{S}$ if $1<\sigma<2$ or $\sigma>2$, respectively. %
Notice that in the limit case of $\sigma=2$ both expressions of $\gamma$ in (\ref{gamma-sigma}) converge to the same value $\gamma=2$. %
Let us set now
\begin{equation}\label{n.n.o.d.e}
\frac{d}{d\,t}E(t)+C_1E(t)^{\gamma}= C_2\Leftrightarrow
\frac{d}{d\,t}E(t)=C_2-C_1E(t)^{\gamma}:=\Lambda(t)\,.
\end{equation}
If $\Lambda(t)<0$ or $\Lambda(t)>0$ at some time $t$ (possibly different), then $E(t)$ is decreasing or increasing, respectively, at that time. %
In consequence, the asymptotically stable equilibrium of (\ref{n.n.o.d.e}) is reached when
\begin{equation}\label{E-ast}
\Lambda(t)=0\Leftrightarrow E(t)=\left(\frac{C_2}{C_1}\right)^{\frac{1}{\gamma}}\equiv
C\left(\frac{C_{\textbf{f}}}{2\nu\min(\nu/2,\alpha)}\right)^{\frac{1}{\gamma}}
:=\frac{\mathcal{E}_{\ast}}{2}\equiv E_{\ast}\,,
\end{equation}
where $C_1$, $C_2$ and $C$ are given in (\ref{C1C2-sigma}). %
A simple analysis of (\ref{n.n.o.d.e}) and (\ref{E-ast}) shows us that if there exists a positive time
$t_0$ such that $E(t_0)<E_{\ast}$, then $0\leq E(t)<E_{\ast}$ for all time $t>t_0$ and, consequently, $E(t)\nearrow E_{\ast}$ as $t\to\infty$. %
In this case, we are done and we obtain for all $t>t_0$
$$\|\textbf{u}(t)\|_{\textbf{L}^2(\Omega)}^2\leq 2C_{GNS}^2
\left(\frac{C_{\textbf{f}}}{2\nu\min(\nu/2,\alpha)}\right)^{\frac{4+(2-\sigma)N}{(2-\sigma)N+2\sigma}}
\quad\mbox{if}\quad 1<\sigma<2$$
or
$$\|\textbf{u}(t)\|_{\textbf{L}^2(\Omega)}^2\leq 2C_{S}
\left(\frac{C_{\textbf{f}}}{2\nu\min(\nu/2,\alpha)}\right)^{\frac{\sigma+2}{2\sigma}}
\quad\mbox{if}\quad \sigma>2\,.$$
The reciprocal case shall be studied in the next theorem.

\begin{theorem}[Exponential decay] \label{main-theorem-ed}
Assume $\textbf{u}_0\in\textbf{H}$ and $\textbf{f}\not=\textbf{0}$ satisfies to (\ref{C-f}). %
Let $\textbf{u}$ be a weak solution of the modified NS problem
(\ref{NSE-2})-(\ref{b-c}) in the sense of Definition~\ref{def-w-sol-tp}. %
In addition assume that there exists a positive time $t_0$ such that $\|\textbf{u}(t_0)\|_{\textbf{L}^2(\Omega)}^2>\mathcal{E}_{\ast}$. %
Then there exists a positive constant $C$ such that
\begin{equation}\label{e-decay}
\|\textbf{u}(t)\|_{\textbf{L}^2(\Omega)}^2\leq
\left(\|\textbf{u}(t_0)\|_{\textbf{L}^2(\Omega)}^2-\mathcal{E}_{\ast}\right)e^{-C(t-t_0)}+\mathcal{E}_{\ast}
\quad\mbox{for all $t>t_0$}\,,
\end{equation}
where $\mathcal{E}_{\ast}$ is given in (\ref{E-ast}).
\end{theorem}
PROOF.
\emph{First Step.}
Again a simple analysis of (\ref{n.n.o.d.e}) and (\ref{E-ast}) shows us that if there exists a positive time
$t_0$ such that $E(t_0)>E_{\ast}$, then $E(t) > E_{\ast}$ for all time $t>t_0$ and, consequently, $E(t)\searrow E_{\ast}$ as $t\to\infty$. %
In order to simplify the notations, let us set $E(t)=E$, $\frac{d}{d\,t}E(t)=E'$ and $f(E):=C_1E^{\gamma}$. %
First, we observe that with these notations we have from (\ref{n.n.o.d.e}) and (\ref{E-ast}),
\begin{equation}\label{o.d.e.E-E*}
E'+C_1E^{\gamma}=C_2\quad\Longleftrightarrow\quad
(E-E_{\ast})'+C_1\frac{E^{\gamma}-E_{\ast}^{\gamma}}{E-E_{\ast}}(E-E_{\ast})=0\,.
\end{equation}

\vspace{0.2cm}\noindent
\emph{Second Step.}
Using (\ref{E-ast}), we can prove that
\begin{equation}\label{frac-g-f}
C_1\frac{E^{\gamma}-E_{\ast}^{\gamma}}{E-E_{\ast}}\equiv\frac{f(E)-f(E_{\ast})}{E-E_{\ast}}>f'(E_{\ast})\equiv\gamma C_1^{\frac{1}{\gamma}}C_2^{\frac{\gamma-1}{\gamma}}\quad\mbox{iff}\quad \gamma>1\,.
\end{equation}
According to the expressions of $\gamma$ (see (\ref{gamma-sigma})), we can see that (\ref{frac-g-f}) holds iff $\sigma>2$. Then from (\ref{o.d.e.E-E*}) and (\ref{frac-g-f}) we derive the following linear differential inequality
\begin{equation}\label{l.d.i-s>2}
(E-E_{\ast})'+C(E-E_{\ast})<0\,,\quad C=\gamma\,C_1^{\frac{1}{\gamma}}C_2^{\frac{\gamma-1}{\gamma}}\,.
\end{equation}
Integrating (\ref{l.d.i-s>2}) between $t_0$ and $t>t_0$, we obtain (\ref{e-decay}) where, from (\ref{C1C2-sigma}) and (\ref{l.d.i-s>2}),
\begin{equation}\label{C-e-decay}
C=\frac{2\sigma}{\sigma+2}\frac{2}{C_S}\min\left(\frac{\nu}{2},\alpha\right)^{\frac{\sigma+2}{2\sigma}}
\left(\frac{C_{\textbf{f}}}{2\nu}\right)^{\frac{\sigma-2}{2\sigma}}\,.
\end{equation}

\vspace{0.2cm}\noindent
\emph{Third Step.} For $1<\sigma<2$ and in addition to what we have done in the First Step, let us set $E(t_0)=E_0$. %
We then prove
\begin{equation}\label{frac-g-ff}
C_1\frac{E^{\gamma}-E_{\ast}^{\gamma}}{E-E_{\ast}}>
C_1\frac{E^{\gamma}_0-E_{\ast}^{\gamma}}{E_0-E_{\ast}}\quad \mbox{iff}\quad 0<\gamma<1\,,
\end{equation}
and according to (\ref{gamma-sigma}), we see that this holds iff $1<\sigma<2$. %
Then from (\ref{o.d.e.E-E*}) and (\ref{frac-g-ff}), we derive the following linear differential inequality
\begin{equation}\label{l.d.i-s<2}
(E-E_{\ast})'+C(E-E_{\ast})<0\,,\quad
C=C_1\frac{E^{\gamma}_0-E_{\ast}^{\gamma}}{E_0-E_{\ast}}\,.
\end{equation}
Integrating (\ref{l.d.i-s<2}) between $t_0$ and $t>t_0$ leads us to (\ref{e-decay}), where from (\ref{C1C2-sigma}) and (\ref{l.d.i-s<2}), the constant $C$ is given by
\begin{equation}\label{CC-e-decay}
C=\left(\frac{2}{C_{GNS}^2}\right)^{\frac{(2-\sigma)N+2\sigma}{4+(2-\sigma)N}}\min\left(\frac{\nu}{2},\alpha\right)
\frac{E_0^{\frac{(2-\sigma)N+2\sigma}{4+(2-\sigma)N}}-E_{\ast}^{\frac{(2-\sigma)N+2\sigma}{4+(2-\sigma)N}}}{E_0-E_{\ast}} \,.\square
\end{equation}

To conclude this section, we analyze the results obtained in this section in terms of $\alpha$ and $\nu$. %
Theorem~\ref{main-theorem} is the most restrictive, because if one and only one of the constants $\alpha$ and $\nu$ tend to zero, we obtain,
from (\ref{t-ast}) and (\ref{epsil-0}), $t_{\ast}\to\infty$ and $\epsilon\to 0$ and those results fail. %
But if only $\alpha$ tend to zero, the results corresponding to exponential and power-time decays (Theorems~\ref{main-theorem-od} and~\ref{main-theorem-ed}) still hold. This situation corresponds to consider the classical Navier-Stokes problem.
If $\nu$ tends to zero, the situation gets worse, because we can no longer use Sobolev-type inequalities.

\section{Conclusions}

Throughout this paper we have seen that the introduction of the absorption term $|\textbf{u}|^{\sigma-2}\textbf{u}$
in the momentum equation of the Navier-Stkes problem allows us to obtain different time properties for the weak solutions. %
If $1<\sigma<2$ we have improved the known results for the NS problem and we obtained an extinction in a finite time, whereas for $\sigma\geq 2$ we have obtained exponential or power-time decays. %
The problem lies in how to justify, from the physical point of view, the appearance of this term in the momentum equation. %
We may think of it as expressing somehow a sink that dissipates kinetics energy. %
Or it can correspond to a physical body inside the fluid, or at the boundary, which makes the flow slower when the time is passing by or even stop it in a finite time. %
Other possibility is to consider the absorption term as resulting from a real or fictitious force as already mentioned in (\ref{force-feedback}). %
For instance, in Geophysical flows the fictitious Coriolis force resulting from the Coriolis acceleration is given by $2\mathbf{\Omega\times}\textbf{u}$, where $\mathbf{\Omega}$ is the Earth angular velocity, 
and $\alpha|\textbf{u}|^{\sigma-2}\textbf{u}$ approximates the Coriolis force if $\sigma=2$. In this case, we only obtain exponential and power-time decays. %
From the theoretical point of view, the results of this paper can be extended to a great variety of Fluid Mechanics problems modified by introducing in the momentum equation an absorption like term. %
In Oliveira~\cite{OB-absorption} we consider the Oberbeck-Boussinesq problem modified by the thermo-absorption term $\alpha|\textbf{u}|^{\sigma(\theta)-2}\textbf{u}$, where $\theta$ is the absolute temperature. %
Here  $\sigma$ is a temperature depending function with Lipschitz regularity and such that $\sigma(\theta)>1$ for all $\theta\in\mathbb{R}$. %
Under the assumption that $1<\sigma^{-}\leq\sigma(\theta)\leq \sigma^{+}<2$ for all $\theta\in\mathbb{R}$,
where $\sigma^{-}$ and $\sigma^{+}$ are constants, we are able to establish the same asymptotic stability properties. %
We may also consider the non-homogeneous Navier-Stokes problem modified by an absorption term. %
In this case, we only need to assume that the density function $\rho$ is bounded:
$\frac{1}{C_{\rho}}\leq \rho\,,\ \rho_0 \leq C_{\rho}$,
where $C_{\rho}$ is a positive constant and $\rho_0$ is the initial density. %
These properties extend also for the modified NS problem supplemented with the slip boundary conditions:
$\mathbf{u\cdot\,n}=0\quad\mbox{and}\quad\mathbf{u\cdot\tau}=\beta^{-1}\mathbf{t\cdot\tau}$ on
$\Gamma_T$, where $\mathbf{n}$ and $\mathbf{\tau}$ denote, respectively, unit
normal and tangential vectors to the boundary $\partial\Omega$,
$\textbf{t}=\mathbf{n\cdot T}$ is the stress vector and
$\beta$ is a coefficient with no defined sign (see Antontsev and Oliveira~\cite{RIMS-2007}). %
As far as the behavior in space of the weak solutions of such modified problems, at the moment, we are able to prove analogous properties only for the 2D stationary modified NS problem (see Antontsev and Oliveira~\cite{JMFM-2004}).

\subsection*{Acknowledgements}
This work was partially supported by FEDER and FCT-Plurianual 2008.

\end{document}